  %

  \input amssym


  \font \bbfive		= bbm5
  \font \bbseven	= bbm7
  \font \bbten		= bbm10
  \font \eightbf	= cmbx8
  \font \eighti		= cmmi8 \skewchar \eighti = '177
  \font \eightit	= cmti8
  \font \eightrm	= cmr8
  \font \eightsl	= cmsl8
  \font \eightsy	= cmsy8 \skewchar \eightsy = '60
  \font \eighttt	= cmtt8 \hyphenchar\eighttt = -1

  \font \sixi		= cmmi6 \skewchar \sixi = '177
  \font \sixrm		= cmr6
  \font \sixsy		= cmsy6 \skewchar \sixsy = '60
  \font \tensc		= cmcsc10

  \font \titlefont	= cmbx12
  \scriptfont \bffam	= \bbseven
  \scriptscriptfont \bffam = \bbfive
  \textfont \bffam	= \bbten

  \newskip \ttglue

  \def \eightpoint {\def \rm {\fam 0 \eightrm }\relax
  \textfont 0= \eightrm
  \scriptfont 0 = \sixrm \scriptscriptfont 0 = \fiverm
  \textfont 1 = \eighti
  \scriptfont 1 = \sixi \scriptscriptfont 1 = \fivei
  \textfont 2 = \eightsy
  \scriptfont 2 = \sixsy \scriptscriptfont 2 = \fivesy
  \textfont 3 = \tenex
  \scriptfont 3 = \tenex \scriptscriptfont 3 = \tenex
  \def \it {\fam \itfam \eightit }\relax
  \textfont \itfam = \eightit
  \def \sl {\fam \slfam \eightsl }\relax
  \textfont \slfam = \eightsl
  \def \bf {\fam \bffam \eightbf }\relax
  \textfont \bffam = \bbseven
  \scriptfont \bffam = \bbfive
  \scriptscriptfont \bffam = \bbfive
  \def \tt {\fam \ttfam \eighttt }\relax
  \textfont \ttfam = \eighttt
  \tt \ttglue = .5em plus.25em minus.15em
  \normalbaselineskip = 9pt
  \def \MF {{\manual opqr}\-{\manual stuq}}\relax
  \let \sc = \sixrm
  \let \big = \eightbig
  \setbox \strutbox = \hbox {\vrule height7pt depth2pt width0pt}\relax
  \normalbaselines \rm }


  \def \ifundef #1{\expandafter \ifx \csname #1\endcsname \relax }


  \newcount \secno \secno = 0
  \newcount \stno \stno = 0
  \newcount \eqcntr \eqcntr = 0

  \def \track #1#2#3{\ifundef {#1}\else \hbox {\sixrm [#2\string #3] }\fi }

  \def \advseqnumbering {\global \advance \stno by 1 \global \eqcntr =0}

  \def \current {\number \secno \ifnum \number \stno = 0 \else .\number \stno \fi }

  \def \laberr #1{\message {*** RELABEL CHECKED FALSE for #1 ***}
      RELABEL CHECKED FALSE FOR #1, EXITING.
      \end }

  \def \deflabel #1#2{%
    \ifundef {#1}%
      \global \expandafter
      \edef \csname #1\endcsname {#2}%
    \else
      \edef \deflabelaux {\expandafter \csname #1\endcsname }%
      \edef \deflabelbux {#2}%
      \ifx \deflabelaux \deflabelbux \else \laberr {#1=(\deflabelaux )=(\deflabelbux )} \fi
      \fi
    \track {showlabel}{*}{#1}}

  \def \eqmark #1 {\advseqnumbering
    \eqno {(\current )}
    \deflabel {#1}{\current }}

  \def \subeqmark #1 {\global \advance \eqcntr by 1
    \edef \subeqmarkaux {\current .\number \eqcntr }
    \eqno {(\subeqmarkaux )}
    \deflabel {#1}{\subeqmarkaux }}

  \def \label #1 {\deflabel {#1}{\current }}
  \def \lcite #1{(#1\track {showlcit}{$\bullet $}{#1})}
  \def \forwardcite #1#2{\deflabel {#1}{#2}\lcite {#2}}


  \catcode `\@=11
  \def \c@itrk #1{{\bf #1}\track {showcitations}{\#}{#1}} 
  \def \c@ite #1{{\rm [\c@itrk{#1}]}}
  \def \sc@ite [#1]#2{[\c@itrk{#2}\hskip 0.7pt:\hskip 2pt #1]}
  \def \du@lcite {\if \pe@k [\expandafter \sc@ite \else \expandafter \c@ite \fi }
  \def \cite {\futurelet \pe@k \du@lcite }
  \catcode `\@=12


  \newcount \bibno \bibno = 0
  \newcount \bibtype \bibtype = 0 


  \def \newbib #1#2{\ifcase \number \bibtype
	\global \advance \bibno by 1 \edef #1{\number \bibno }\or
	\edef #1{#2}\or
	\edef #1{\string #1}\fi }

  \def \bibitem #1#2#3#4{\smallbreak \item {[#1]} #2, ``#3'', #4.}

  \def \references {\begingroup \bigbreak \eightpoint
    \centerline {\tensc References}
    \nobreak \medskip \frenchspacing }


  \def \Headlines #1#2{\nopagenumbers
    \headline {\ifnum \pageno = 1 \hfil
    \else \ifodd \pageno \tensc \hfil \lcase {#1} \hfil \folio
    \else \tensc \folio \hfil \lcase {#2} \hfil
    \fi \fi }}

  \long \def \Quote #1\endQuote {\begingroup \leftskip 35pt \rightskip 35pt
\parindent 17pt \eightpoint #1\par \endgroup }
  \long \def \Abstract #1\endAbstract {\bigskip \Quote \noindent #1\endQuote }

  \def \Note #1{\footnote {}{\eightpoint #1}}
  \def \Date #1 {\Note {\it Date: #1.}}


  \def \lcase #1{\edef \auxvar {\lowercase {#1}}\auxvar }

  \def \section #1 \par {\global \advance \secno by 1 \stno = 0
    \bigbreak \noindent {\bf \number \secno .\enspace #1.}
    \nobreak \medskip \noindent }

  \def \state #1 #2\par {\medbreak \noindent \advseqnumbering {\bf \current .\enspace #1.\enspace \sl #2\par }\medbreak }
  \def \definition #1\par {\state Definition \rm #1\par }

  \long \def \Proof #1\endProof {\medbreak \noindent {\it Proof.\enspace }#1
\ifmmode \eqno \endproofmarker $$ \else \hfill $\endproofmarker $ \looseness = -1 \fi \medbreak }

  \def \$#1{#1 $$$$ #1}
  \def \explain #1#2{\mathrel {\buildrel \hbox {\sixrm (#1)} \over #2}}
  \def \=#1{\explain {#1}{=}}

  \def \pilar #1{\vrule height #1 width 0pt}

  \newcount \footno \footno = 1
  \newcount \halffootno \footno = 1
  \def \footcntr {\global \advance \footno by 1
  \halffootno =\footno
  \divide \halffootno by 2
  $^{\number \halffootno }$}
  \def \fn #1{\footnote {\footcntr }{\eightpoint #1\par }}




  \def \Item #1{\smallskip \item {{\rm #1}}}
  \newcount \zitemno \zitemno = 0

  \def \izitem {\global \zitemno = 0}
  \def \zitemplus {\global \advance \zitemno by 1 \relax }
  \def \rzitem {\romannumeral \zitemno }
  \def \rzitemplus {\zitemplus \rzitem } 
  \def \zitem {\Item {{\rm (\rzitemplus )}}}
  \def \Zitem {\izitem \zitem }
  \def \zitemmark #1 {\deflabel {#1}{\rzitem }}

  \newcount \nitemno \nitemno = 0
  
  \def \nitem {\global \advance \nitemno by 1 \Item {{\rm (\number \nitemno )}}}

  \newcount \aitemno \aitemno = -1
  \def \boxlet #1{\hbox to 6.5pt{\hfill #1\hfill }}
  \def \iaitem {\aitemno = -1}
  \def \aitemconv {\ifcase \aitemno a\or b\or c\or d\or e\or f\or g\or
h\or i\or j\or k\or l\or m\or n\or o\or p\or q\or r\or s\or t\or u\or
v\or w\or x\or y\or z\else zzz\fi }
  \def \aitem {\global \advance \aitemno by 1\Item {(\boxlet \aitemconv )}}
  \def \aitemmark #1 {\deflabel {#1}{\aitemconv }}


  \font \mf =cmex10
  \def \union {\mathop {\raise 9pt \hbox {\mf S}}\limits }
  \def \inters {\mathop {\raise 9pt \hbox {\mf T}}\limits }

  \def \<{\left \langle \vrule width 0pt depth 0pt height 8pt }
  \def \>{\right \rangle }
  \def \({\big (}
  \def \){\big )}
  
  \def \and {\hbox {,\quad and \quad }}

  \def \imply {\kern 7pt \Rightarrow \kern 7pt}
  \def \for #1{,\quad \forall \,#1}
  \def \endproofmarker {\square } 
  \def \"#1{{\it #1}\/} 
  
  \def \*{\otimes }
  \def \caldef #1{\global \expandafter \edef \csname #1\endcsname {{\cal #1}}}
  \def \bfdef #1{\global \expandafter \edef \csname #1\endcsname {{\bf #1}}}
  \bfdef N \bfdef Z \bfdef C \bfdef R


  \def \Caixa #1{\setbox 1=\hbox {$#1$\kern 1pt}\global \edef \tamcaixa {\the \wd 1}\box 1}
  \def \caixa #1{\hbox to \tamcaixa {$#1$\hfil }}



  \catcode `\@=11

  \def \overparenOnefill {$\m@th
  \setbox 0=\hbox {$\braceld $}%
  \braceld \leaders \vrule height\ht 0 depth0pt\hfill
  \leaders \vrule height\ht 0 depth0pt\hfill \bracerd $}

  \def \overparenOne #1{\mathop {\vbox {\m@th\ialign {##\crcr \noalign {\kern -1pt}
  \overparenOnefill \crcr \noalign {\kern 3pt\nointerlineskip }
  $\hfil \displaystyle {#1}\hfil $\crcr }}}\limits }

  \def \overparenTwofill {$\m@th
  \lower 0.3pt \hbox {$\braceld $}
  \leaders \vrule depth 0pt height1pt \hfill
  \lower 0.3pt \hbox {$\bracerd $}$}

  \def \overparenTwo #1{\mathop {\vbox {\ialign {##\crcr \noalign {\kern -1pt}
  \overparenTwofill \crcr \noalign {\kern 3pt\nointerlineskip }
  $\hfil \displaystyle {#1}\hfil $\crcr }}}\limits }

  \catcode `\@=12

  \def\soma{\mathop {\textstyle\sum}\limits}
  \def\su#1#2{\soma_{#1=1}^#2}
  \def\sumn{\su in}
  \def\summ{\su jm}

  \def\frac#1#2{{\raise -2pt \hbox{$\scriptstyle #1$}\over {\raise 2pt \hbox{$\scriptstyle #2$}}}}
  \def\frac#1#2{
    \edef\abaixa{\ifinner 0pt \else 2pt \fi}
    {\raise -\abaixa \hbox{$\scriptstyle #1$}\over {\raise \abaixa \hbox{$\scriptstyle #2$}}}}

  \def\rtinv#1{{\scriptstyle  1\over \scriptstyle \sqrt#1}}
  \def\rtinv#1{\frac 1{\sqrt#1}}

  \def\A{A}
  \def\Ap{\A_p}
  \def\hAp{\hat\A_p}
  \def\Aq{\A_q}

  \catcode`\@=11
  \def\narrowmatrix#1{\null\,\vcenter{\normalbaselines\m@th
	\ialign{\hfil$##$\hfil&&\ \hfil$##$\hfil\crcr
	\mathstrut\crcr\noalign{\kern-\baselineskip}
	#1\crcr\mathstrut\crcr\noalign{\kern-\baselineskip}}}\,}
  \catcode`\@=12

  \def\bool#1{{\scriptstyle [#1]}\,}
  \def\Lnm{L_{n,m}}
  \def\Onm{{\cal O}_{n,m}}
  \def\Onmp{\Onm^p}   \def\Onmq{\Onm^q}
  \caldef V \caldef H \caldef G \caldef F
  \def\S{{\cal S}_{n,m}}
  \def\R{{\cal R}_{n,m}}
  \def\X{{\cal X}}
  \def\u#1{\hat{#1}}
  \def\uo#1{\check{#1}}
    
  \def\pio{\uo\pi}
  \def\So{\uo  s}
  \def\piu{\u\pi}
  \def\Su{\u s}
  \caldef T
  \def\Top#1{\T(#1,\V,\H)}
  \def\CovAlg#1{C^*\(#1,\V,\H\)}    \def\CovAlg#1{#1\kern-2pt \rtimes_{_{\V,\H}}\kern-2pt\N}
  \def\KV{{\cal K}_\V}
  \def\KH{{\cal K}_\H}
  \def\resp#1{{\rm[resp.}~#1{\rm]}}
  \def\leftright{left- \resp{right-}}
  \def\mat#1#2#3#4{\pmatrix{#1 & #2 \cr #3 & #4}}
  \def\lmat#1{\mat00{#1}0}
  
  \def\clspan{\overline{\rm span}\ }
  \def\fix{\smallskip\noindent $\blacktriangleright$\kern 12pt }
  \def\hp{\hat p} \def\hq{\hat q}



\bibtype 0

\newbib {\AraExel}{AEK}
\newbib {\Arz}{A}
\newbib {\ArzRena}{AR}
\newbib {\RieffelMorita}{BGR}
\newbib {\CuVer}{CV}
\newbib {\tpa}{E1}
\newbib {\amena}{E2}
\newbib {\inverse}{E3}
\newbib {\ortho}{E4}
\newbib {\ExelEndo}{E5}
\newbib {\interactions}{E6}
\newbib {\newsgrp}{E7}
\newbib {\oberwolfach}{E8}
\newbib {\infinoa}{EL}
\newbib {\poly}{ER}
\newbib {\Zettl}{Z}

  \font\titlefont=cmbx10
  \font\titlefontsl=cmbxsl10

  \Headlines {INTERACTIONS AND DYNAMICAL SYSTEMS OF TYPE $(n,m)$}
    {R.~Exel}
  \null\vskip -1cm
  \centerline{\titlefont INTERACTIONS AND DYNAMICAL SYSTEMS}
  \smallskip
  \centerline{\titlefont OF TYPE {\titlefontsl (n,m)} -- A CASE STUDY}

  \Date{22 December 2012}

  \footnote{\null}
  {\eightrm 2010 \eightsl Mathematics Subject Classification:
  \eightrm
  46L05, 
  46L55. 
  }

  \Note
  {\it Key words and phrases: \rm Leavitt C*-algebra, $\Onm$, interactions.}

  \bigskip
  \centerline{\tensc Ruy Exel}

  \Note{Partially supported by CNPq.}

  \bigskip

\Abstract We prove that the C*-algebra of the universal $(n,m)$-dynamical system may be obtained, up to Morita--Rieffel
equivalence, as the crossed-product relative to an interaction on a commutative C*-algebra.  The interaction involved is
shown not to be part of an interaction group.  \endAbstract

\section Introduction

  The notion of \"{interactions} was introduced in \cite{\interactions} in order to provide a common generalization for
endomorphisms of C*-algebras and their transfer operators.  One of the main results in \cite{\interactions}, namely
Theorem 6.3, is the proof of the existence of a \"{covariant representation} for any given interaction, but no
consistent notion of \"{crossed product} was introduced.

In reality, in the last section of \cite{\interactions}, an admittedly experimental attempt was made to provide some
sort of crossed product in terms of a certain generalization of the Cuntz--Pimsner algebra to a context in which the
correspondence is replaced by a \"{generalized correspondence} \cite[Definition 7.1]{\interactions}.  However, no
nontrivial examples were provided so the theory was not put through any significant test.

The notion of interactions was later given a (non-equivalent) alternative form in \cite{\newsgrp} (see also
\cite{\oberwolfach}), the catch-word being \"{interaction groups}, and a well developed notion of crossed product was
introduced.  Several examples were later exhibited in \cite{\poly}, including the case of the \"{multi-valued map}
$z\mapsto z^{2/3}$ on the circle, which has not received a lot of attention in the literature, except for the paper
\cite{\ArzRena} by Arzumanian and Renault (which was in fact slightly rectified by \cite{\poly}) and some recent work by
Arzumanian \cite{\Arz} and by Cuntz and Vershik \cite{\CuVer}.

In a very rough sense, interactions are related to partial isometries, while interaction groups are related to power
partial isometries, meaning partial isometries whose powers are still partial isometries.  Since partial isometries not
satisfying the latter property are hard to study objects, I eventually developed the impression that interactions
should be likewise considered.

Roughly five years after the appearance of \cite{\interactions}, I was involved in a seemingly unrelated joint project
with P.~Ara and T.~Katsura \cite{\AraExel}, where we introduced the notion of \"{$(n,m)$-dynamical systems} and their
accompanying C*-algebras, denoted $\Onm$, which turned out to be a generalization of the Cuntz algebras.  The method
used to study $\Onm$ was based on partial actions and in no moment did it occur to us to study it from the point of view
of interactions.

By an $(n,m)$-dynamical system we mean two compact spaces $X$ and $Y$, with maps
  $$
  h_1,\ldots,h_n,\ v_1,\ldots,v_m: Y \to X,
  $$
  which are homeomorphisms onto their ranges, and such that
  $$
  X = \bigcup_{i=1}^n h_i(Y) = \bigcup_{j=1}^n v_j(Y),
  $$
  both unions being disjoint unions.  Given such a system, one may consider a local homeomorphism
  $
  \alpha : X \to Y,
  $
  defined to coincide with $h_i^{-1}$ on the range of each $h_i$.  One might consider $\alpha$ as some version of Bernoulli's
shift, for which the $h_i$'s are the inverse branches.

Replacing the $h_i$ by the $v_j$, one may similarly define another local homeomorphism, say $\beta: X \to Y$, having the $v_j$
as inverse branches.

Evidently neither $\alpha$ nor $\beta$ are invertible (unless $n$ or $m=1$), but we may view the \"{multi-valued} map
  $$
  L : y \mapsto \{h_1(y),\ldots,h_n(y)\},
  $$
  as playing the role of the inverse of $\alpha$.  Likewise
  $$
  M : y \mapsto \{v_1(y),\ldots,v_m(y)\}
  $$
  may be considered as some sort of inverse for $\beta$. Playing in a totally careless way with these maps, one may define
  $$
  \V,\H:X\to X,
  $$
  by $\V=M\alpha$, and $\H=L\beta$, and argue that
  $$
  \V^{-1} = \alpha^{-1}M^{-1} = L\beta = \H.
  $$

  Evidently all of this is nonsense, but the notion of interactions may give it a precise and meaningful treatment.  The
main idea is that, when a map is multivalued, it defines a (singly valued) map on the algebra of continuous functions by
averaging over the multiple values.  In the present case, this leads us to defining maps $\V$ and $\H$ on the algebra
$C(X)$ by
  $$
  \V(f)|_x = \frac 1m\summ f\big(v_j(\alpha(x)\big)
  \and
  \H(f)|_x = \frac 1n\sumn f\big(h_j(\beta(x)\big),
  \eqno{(\dagger)}
  $$
  for all $f\in C(X)$, and all $x\in X$.

The interesting fact is that the pair $(\V,\H)$ turns out to be an interaction and, moreover, the experimental notion of
crossed product introduced in \cite{\interactions} fits like a glove in the present situation, producing the expected
result, namely the full hereditary subalgebra of $\Onm$ associated to the characteristic function on $X$.

Besides briefly recalling the necessary background, the content of this paper is precisely to prove the isomorphism of
the crossed product $\CovAlg{C(X)}$ with the hereditary subalgebra of $\Onm$ mentioned above.

Another question that we discuss is the possibility of fitting the theory of interaction groups to $\Onm$, but we
unfortunately find in \forwardcite{NotPower}{3.11} that this is not possible.

Before we actually begin, we should say that the description of $\V$ and $\H$ given in $(\dagger)$, above, is not quite the
one we use below, as we have chosen to emphasize the algebraic aspects of $\Onm$ over its dynamical picture.  However,
without too much effort, the reader may use the results in \cite{\AraExel} to show that $(\dagger)$ agrees with the
definitions of $\V$ and $\H$ given in \forwardcite{introVH}{3.7}, below.


\section Interactions

In this section we will give a brief overview of the notions of interactions and the corresponding crossed product.  For
more information the reader is referred to \cite{\interactions}.

\fix From now on we will let $A$ be a fixed unital C*-algebra.

  \definition \cite[Definition 3.1]{\interactions}\quad \label DefineInteraction
  A pair $(\V,\H)$ of maps
  $$
  \V,\H:A\to A
  $$
  will be called an \"{interaction over $A$}, if
  \izitem
  \zitem $\V$ and $\H$ are positive, bounded, unital linear maps,
  \zitem $\V\H\V=\V$,
  \zitem $\H\V\H=\H$,
  \zitem $\V(xy)=\V(x)\V(y)$, if either $x$ or $y$ belong to $\H(A)$,
  \zitem $\H(xy)=\H(x)\H(y)$, if either $x$ or $y$ belong to $\V(A)$.

\fix Let us assume, for the remainder of this section, that $(\V,\H)$ is a fixed interaction over $A$.

\definition \cite[Definition 3.5]{\interactions}\quad A
  \"{covariant representation} of $(\V,\H)$ in a given unital C*-algebra $B$ is a pair $(\pi, s)$, where $\pi$ is a unital
*-homomorphism of $A$ into $B$, and $s$ is a partial isometry in $B$ such that
  \izitem
  \zitem $s\pi(a)s^* = \pi(\V(a)) ss^*$, and
  \zitem $s^*\pi(a) s = \pi(\H(a)) s^*s$,
  \medskip\noindent for every $a$ in $A$.

\definition We will denote by $\Top A$ the universal
  unital\fn{When we say \"{universal unital} we mean that we are working in the category of unital C*-algebras and hence
all algebras and morphisms involved in its universal properties are supposed to be unital.}
  C*-algebra generated by a copy of $A$ and a partial isometry $\So$, subject to the relations
  \izitem
  \zitem $\So a\So^* = \V(a)\So\So^*$, and
  \zitem $\So^*a\So = \H(a)\So^*\So$,
  \medskip\noindent for every $a$ in $A$.
  The canonical mapping from $A$ to $\Top A$ will be denoted by $\pio$.

It is readily seen that $\(\pio,\So \)$ is a covariant representation of $(\V,\H)$ in $\Top A$.  In addition, $\Top A$
is clearly the universal C*-algebra for covariant representations of $(\V,\H)$ in the sense that any covariant
representation factors through $\Top A$.

We should remark that, as we are working in the category of unital C*-algebras and morphisms, the natural inclusion
$\pio$ of $A$ in $\Top A$ is a unital map and, in particular,
  $$
  \pio(1)\So = \So\pio(1) = \So.
  \eqmark UnitalTop
  $$

  Quite likely is is also possible to develop a similar theory for non-unital algebras but, given the examples we have
in mind, we have decided to concentrate on the unital case here.

\state Proposition \cite{\interactions}\quad The closed linear span of $\pio(A)\So\pio(A)$, henceforth denoted by $\X$,
is a ternary ring of operators \cite{\Zettl}, meaning that it satisfies
  $$
  \X\X^*\X\subseteq\X.
  $$

  \Proof
  For all $a,b,c,d,e,f\in A$, we have
  $$
  \(\pio(a)\So \pio(b)\)\(\pio(c)\So \pio(d)\)^*\(\pio(e)\So \pio(f)\) =
  \pio(a)\So \pio(bd^*)\So^*\pio(c^* e)\So \pio(f) \$=
  \pio(a) \pio\(\V(bd^*)\)\So^*\So\pio\(\H(c^* e)\) \pio(f) =
  \pio\(a\V(bd^*)\)\So\pio\(\H(c^* e)f\) \in\X.
  \endProof

From the above result it follows that
  $$
  \KV := \clspan \X\X^*
  $$
  as well as
  $$
  \KH := \clspan \X^*\X
  $$
  are closed *-subalgebras of $\Top A$.  It also follows that
  $$
  \KV\X\subseteq\X \and \X\KH \subseteq\X,
  $$
  and hence that $\X$ is a $\KV$--$\KH$--bimodule.  On the other hand it is easily seen that $\X$ is an $A$--$A$--bimodule.

\definition \cite{\interactions}
  \iaitem
  \aitem A
  \"{{\leftright} redundancy} is a pair
  $(a,k)$ in $A\times\KV$ \resp{$A\times\KH$}
  such that
  $$
  \pio(a)x = kx\quad \resp{x\pio(a) = xk} \for x \in\X.
  $$
  \aitem
  The \"{redundancy ideal} is the closed two-sided ideal
  of $A$ generated by the set
  $$
  \big\{\pio(a)-k: (a,k) \hbox{ is a left-redundancy}\big\} \cup
  \big\{\pio(a)-k: (a,k) \hbox{ is a right-redundancy}\big\}.
  $$
  \aitem
  The quotient of $\Top A$ by the redundancy ideal will be called the \"{covariance algebra} or the \"{crossed product}
for the interaction $(\V,\H)$, and will be denoted by $\CovAlg A$.
  \aitem Letting
  $$
  q: \Top A \to\CovAlg A
  $$
  be the quotient map, we will let $\piu = q\circ\pio$, and $\Su = q(\So)$.

\bigskip Again we have that $(\piu,\Su)$ is a covariant representation of $(\V,\H)$ in $\CovAlg A$.

The following is an elementary result which slightly simplifies some computations involving redundancies:

\state Proposition \label SimpleRedundancies
  A pair $(a,k)$ in $A\times\KV$ \resp{$A\times\KH$} is a {\leftright} redundancy if and only if
  $$
  \pio(a)\pio(b)\So = k\pio(b)\So\quad \resp{\So \pio(b)\pio(a) = \So \pio(b)k} \for a \in A.
  $$

\Proof This follows immediately from \lcite{\UnitalTop} and the density of $\pio(A)\So\pio(A)$ in $\X$. \endProof

\section Brief description of $\Onm$

Let us now introduce the algebra $\Onm$ which will play a prominent role in our main result below.  For further details
on the properties and structure of $\Onm$ the reader is referred to \cite{\AraExel}.

\definition \label DefineOnm Given integers $n,m\geq1$, the \"{Leavitt C*-algebra of type $(n,m)$}, henceforth denoted by
$\Lnm$, is the universal unital C*-algebra generated by partial isometries
  $s_1,\ldots,s_n, \ t_1,\ldots,t_m$
  satisfying the relations
  $$
  s_i^*s_k = 0, \hbox{ for } i\neq k, $$$$
  t_j^*t_l = 0, \hbox{ for } j\neq l, $$$$
  s_i^*s_i = t_j^*t_j =: q, $$$$
  \sum_{i=1}^n s_is_i^* = \sum_{j=1}^m t_jt_j^* =:p, $$$$
  pq=0, \quad p+q=1.
  $$

As observed in \cite[Section 2]{\AraExel}, when $n,m>1$, the partial isometries $s_i$ and $t_j$ in $\Lnm$ do not form a
\"{tame} set, in the sense that the multiplicative subsemigroup of $\Lnm$ generated by
  $$
  \big\{s_1, \ldots, s_n,\ s_1^*, \ldots, s^*_n,\ t_1, \ldots, t_m,\ t_1^*, \ldots, t^*_m\big\}
  $$
  does not consist of partial isometries.  In order to fix this, we consider the ideal $J\trianglelefteq\Lnm$ generated
by all elements of the form $xx^*x-x$, where $x$ runs in the above mentioned semigroup.

\definition $\Onm$ is the quotient of $\Lnm$ by the ideal $J$ described above.

From now on we will concentrate our attention on $\Onm$, whereas $\Lnm$ will not play any further role in this work.  We
will therefore not bother to introduce any new notation for the images of the $s_i$ and $t_j$ in $\Onm$, denoting them
again by $s_i$ and $t_j$, as no confusion will arise.

\definition \label IntroS The multiplicative subsemigroup of $\Onm$ generated by all the $s_i$, all the $t_j$, as well
as their adjoints, will be denoted by $\S$.

It is then clear that $\S$ is formed by partial isometries, and hence it is an \"{inverse semigroup}.  Its \"{idempotent
semi-lattice}, namely
  $$
  E(\S) = \{ss^*: s\in\S\} = \{e\in\S: e^2=e\}
  $$
  is a set of commuting projections and therefore generates an abelian sub-C*-algebra of $\Onm$, which we will denote by
$\A$.  Sections (2) and (4) of \cite{\AraExel} give two different descriptions of the spectrum of $\A$.

Evidently $p$ and $q$ are complementary (central) projections in $\A$, so $\A$ admits a decomposition as a direct sum of
two ideals:
  $$
  \A = \Ap \oplus \Aq,
  $$
  where $\Ap = p\A$, and $\Aq =q\A$.  Observing that
  $$
  s_i=ps_iq \and t_j = pt_jq,
  \eqmark sEqualpsq
  $$
  for all $i$ and $j$, we see that the ideals generated by either $p$ or $q$ in $\Onm$ coincide with the whole of
$\Onm$, which is to say that $p$ and $q$ are \"{full} projections.  Consequently the subalgebras of $\Onm$ given by
  $$
  \Onmp := p(\Onm) p \and \Onmq := q(\Onm) q
  $$
  are both \"{full corners}, and hence Morita--Rieffel equivalent \cite[Theorem 1.1]{\RieffelMorita} to $\Onm$.

\state Proposition
  For every $i\leq n$, and $j\leq m$, the correspondences
  $$
  \alpha_i: g \mapsto s_i g s_i^*
  \and
  \beta_j: g \mapsto t_j g t_j^*
  $$
  give well defined *-homomorphisms from $\Aq$ to $\Ap$, and the same is true with respect to
  $$
  \alpha := \sumn \alpha_i
  \and
  \beta := \summ \beta_j.
  $$
  Moreover $\alpha$ and $\beta$ are unital.

  \Proof Left to the reader. \endProof

The *-homomorphisms $\alpha_i$ and $\beta_j$ above are closely related to a \"{partial action} of the free group ${\bf F}_{n+m}$
on $\A$ which contains enough information to reconstruct $\Onm$ in the sense that $\Onm$ is isomorphic to the crossed
product $\A\ifundef {rtimes} \times \else \rtimes \fi{\bf F}_{n+m}$.  See \cite{\AraExel} for more information on this.

Another easy consequence of the relations defining $\Onm$ above is in order.

\state Proposition \label IntroML Define maps $L,M:\Ap \to \Aq$ by
  $$
  L(f) = \frac 1n \sumn s_i^* f s_i
  \and
  M(f) = \frac 1m \summ t_j^* f t_j.
  $$
  Then $L$ and $M$ are unital positive linear maps and moreover
  \Zitem $L\alpha$ and $M\beta$ coincide with the identity of $\Aq$.
  \zitem $L\(\alpha(g) f\) = g L(f)$, for all $g \in\Aq$, and $f\in\Ap $.
  \zitem $M\(\beta(g) f\) = g M(f)$, for all $g \in\Aq$, and $f\in\Ap $.

\Proof
  It is clear that $L$ and $M$ are positive linear maps.  Observing that $p$ is the unit of $\Ap$, we have that
  $$
  L(p) = \frac 1n \sumn \su jn s_i^*s_js_j^*s_i = \frac 1n \sumn s_i^*s_i = q,
  $$
  which is the unit of $\Aq$.  Therefore $L$ is indeed a unital map, and a similar argument applies to prove that $M$ is
also unital.
  In order to prove (ii), let $ g \in\Aq$, and $ f\in\Ap$.  Then
  $$
  L\(\alpha(g) f\) =
  \frac 1n \sumn \su jn s_i^*s_j g s_j^* f s_i =
  \frac 1n \sumn s_i^*s_i g s_i^* f s_i =
  \frac 1n \sumn g s_i^* f s_i =
  g L(f).
  $$
  The proof of (iii) follows along similar lines and, finally, (i) follows from (ii) and (iii) upon plugging $ f=1$.
\endProof

Notice that equations \lcite{\IntroML.ii--iii} bear a close similarity with the axioms defining transfer operators in
\cite{\ExelEndo}.

\state Proposition \label introVH Let
  $\V$ and $\H$ be the linear operators on $\Ap$ defined by
  $$
  \V=\alpha M \and \H=\beta L.
  $$
  Then $(\V,\H)$ is an interaction over $\Ap$.


  \Proof
  It is clear that $\V$ and $\H$ are bounded positive linear maps. In order to prove \lcite{\DefineInteraction.ii}, we
have by \lcite{\IntroML} that
  $$
  \V\H\V = \alpha M\beta L\alpha M = \alpha M = \V.
  $$
  The proof of \lcite{\DefineInteraction.iii} is similar.  As for \lcite{\DefineInteraction.v}, let $ f_1, f_2\in\Ap $, with
$f_1\in\V(\Ap)$.  Then there is $k\in\Ap $ such that
  $$
  f_1 = \V(k) = \alpha\(M(k)\) = \alpha(g),
  $$
  where $g = M(k) \in\Aq$.  We then have by \lcite{\IntroML} that
  $$
  \H(f_1 f_2) = \beta\(L\(\alpha(g) f_2\)\) = \beta\(gL(f_2)\) = \beta(g) \beta \(L(f_2)\) = \beta\(M(k)\) \H(f_2) = \cdots
  $$
  Noticing that
  $
  \H(f_1) = \beta L \alpha M(k) = \beta\(M(k)\),
  $
  the above equals
  $$
  \cdots= \H(f_1)\H(f_2),
  $$
  proving \lcite{\DefineInteraction.v}.
  The proof of \lcite{\DefineInteraction.iv} is similar.
  \endProof

Our next goal will be to produce a covariant representation of $(\V, \H)$ in $\Onmp$.  The partial isometry involved
will actually be produced in terms of two other partial isometries, as follows:

\state Proposition \label IntroPisoR Let
  $$
  S = \rtinv n \sumn s_i
  \and
  T = \rtinv m \summ t_j.
  $$
  Then $S^*S = T^*T = q$.  Consequently $S$ and $T$ are partial isometries.  In addition
  $$
  R := ST^*
  $$
  is a partial isometry belonging to $\Onmp$, which satisfies
  $RR^* = SS^*$ and $R^*R = TT^*$.

  \Proof
  We have
  $$
  S^*S =
  \frac 1n \sumn \su jn s_i^*s_j =
  \frac 1n \sumn s_i^*s_i = q,
  $$
  and similarly $T^*T = q$.  This shows that $S$ and $T$ are partial isometries.  In order to show that $R$ is also a
partial isometry we compute
  $$
  RR^*R = ST^* TS^* ST^* = SqT^* = ST^* = R.
  $$
  By \lcite{\sEqualpsq} we have that $S = pS$, and $T = pT$, so
  $$
  R = ST^* = pST^*p \in \Onmp.
  \endProof

The partial isometries $S$ and $T$ above have a close relationship with the maps $\alpha$, $\beta$, $L$ and $M$ studied above, as
we shall now see.

\state Proposition \label SLMab
  For every $ f\in\Ap$, and for every $ g \in\Aq$, one has that
  \izitem
  \zitem $S^* f S = L(f)$,
  \zitem $T^* f T = M(f)$,
  \zitem $S g = \alpha(g)S $,
  \zitem $T g = \beta(g)T $.

\Proof For $f\in\Ap$, and $i\neq j$, one has that
  $$
  s_i^* f s_j =
  s_i^*s_is_i^* f s_j =
  s_i^* f s_is_i^* s_j = 0.
  $$
  Therefore
  $$
  S^* f S =
  \frac 1n\sumn\su jn s_i^* f s_j =
  \frac 1n\sumn s_i^* f s_i =
  L(f),
  $$
  proving (i).  A similar argument proves (ii).  Given $g \in\Aq$, we have that
  $$
  \alpha(g)S =
  \sumn s_i g s_i^* \rtinv n \su jn s_j =
  \rtinv n \sumn \su jn s_i g s_i^*s_j =
  \rtinv n \sumn s_i g q =
  S g,
  $$
  proving (iii).  A similar computation proves (iv).
  \endProof

The similarity of \lcite{\SLMab.i--iv} with the axioms defining covariant representations in the context of endomorphisms
and transfer operators \cite{\ExelEndo} should again be noticed.

The covariant representation announced above may now be presented.

\state Proposition \label MapOnTop Let $\iota$ denote the inclusion of $\Ap $ into $\Onmp$.  Then $(\iota,R)$ is a covariant
representation of the interaction $(\V,\H)$ in $\Onmp$.  Therefore there is a *-homomorphism
  $$
  \Phi : \Top \Ap \to \Onmp
  $$
  satisfying $\Phi\big(\pio(a)\big) = a$, for all $a$ in $\Ap$, and such that $\Phi(\So) = R$.

\Proof
  Given $f \in \Ap$, we have by \lcite{\SLMab} that
  $$
  Rf R^* = ST^*fTS^* = SM(f)S^* = \alpha\(M(f)\)SS^* = \V(f) RR^*.
  $$
  while
  $$
  R^* f R = TS^*fST^* = TL(f)T^* = \beta\(L(f)\)TT^* = \H(f) R^*R.
  $$
  This shows that $(\iota,R)$ is indeed a covariant representation.
  The last sentence of the statement now follows immediately from the universal property of $\Top \Ap$.
  \endProof

Let us make a short pause to compare the representation above with the representations arising from the theory of
interaction groups \cite{\newsgrp}.  Observe that the partial isometries $v_g$ of \cite[Definition 4.1]{\newsgrp} lie in
the range of a *-partial representation.  Moreover, by \cite[Proposition 2.4.iii]{\inverse}, any two partial isometries
belonging to the range of the same partial representation must have commuting range projections.  In particular, every
such partial isometry is necessarily a \"{power partial isometry}, meaning that its powers are still partial isometries.

\state Proposition \label NotPower If $n$ and $m$ are both greater or equal to 2, the partial isometry $R$ introduced in
\lcite{\IntroPisoR} is not a power partial isometry.  More precisely, $R^2$ is not a partial isometry.

\Proof
  It is well known (see e.g.~\cite[Lemma 5.3]{\ortho}) that the product of two partial isometries $u$ and $v$ is a
partial isometry if and only if the source projection of $u$ commutes with the range projection of $v$.  Thus, $R^2$ is a
partial isometry if and only if $R^*R$ commutes with $RR^*$.  In view of the last sentence of \lcite{\IntroPisoR}, we must
check whether or not $SS^*$ commutes with $TT^*$.
  We have
  $$
  SS^*TT^* = \frac 1{nm}\su{i,k}n\ \su{j,l}m s_is^*_kt_jt^*_l,
  $$
  while
  $$
  TT^*SS^* = \frac 1{nm} \su{j,l}m \ \su{i,k}n t_jt^*_ls_is^*_k.
  $$
  Using the description of $\Onm$ as a partial crossed product \cite[2.5]{\AraExel}, and also the fact that the crossed
product may be defined \cite[Section 2]{\infinoa} as the as the cross sectional C*-algebra of the semi-direct product
Fell bundle \cite[Definition 2.8]{\tpa}, we deduce that $\Onm$ is a cross sectional algebra for a Fell bundle over the
free group ${\bf F}_{n+m}$.

Moreover, if the generators of ${\bf F}_{n+m}$ are denoted $a_1,\ldots,a_n,\ b_1,\ldots,b_m$, each summand $s_is^*_kt_jt^*_l$ in the
expression for $SS^*TT^*$ above lie in the homogeneous space associated to the group element $\pilar{12pt}a_ia^{-1}_kb_jb^{-1}_l$,
and a similar fact holds for the terms $t_jt^*_ls_is^*_k$ in the expression for $TT^*SS^*$.  Should $SS^*$ commute with
$TT^*$, the \"{Fourier coefficient} \cite[Definition 2.7]{\amena} of $SS^*TT^*$ relative to the group element $a_1a_2^{-1}b_1b^{-1}$
would be zero, as this is clearly the case for $TT^*SS^*$.  This means that $s_1s_2^*t_1t_2^*=0$, a contradiction.
  \endProof

As already discussed before the statement of the Proposition above, the fact that $R$ is not a power partial isometry
says that it is impossible to view the covariant representation given by \lcite{\MapOnTop} as part of some interaction
group.

Our next task will be to show that $\Phi$ vanishes on the redundancy ideal of $\Top \Ap$.  The following technical result
initiates our preparations for this.

\state Proposition \label TechLemma Let $y\in\Onmp$.
  \Zitem If \ $y\Ap R \phantom{^*} =\{0\}$, then $y=0$.
  \zitem If \ $y\Ap R^* =\{0\}$, then $y=0$.

\Proof (i) For any given $k\leq n$, notice that $s_ks_k^*\in \Ap$, so
  $$
  0 =
  y s_ks_k^* R =
  y s_ks_k^* ST^* =
  ys_ks_k^*\big(\frac 1{\sqrt{nm}} \sumn\summ s_it^*_j\big) =
  \frac 1{\sqrt{nm}} \sumn\summ y s_ks_k^*s_it^*_j =
  \frac 1{\sqrt{nm}} \summ y s_kt^*_j.
  $$
  Multiplying this on the right by $t_1$, we deduce that
  $$
  0 = \summ y s_kt^*_jt_1 =
  y s_kt_1^*t_1 =
  y s_kq = ys_k.
  $$
  Therefore
  $$
  yp = \su kn ys_ks_k^* = 0.
  $$
  Since $y\in \Onmp$ by hypothesis, we have that $y = yp = 0$.
  The proof of (ii) is similar.
  \endProof

We may now show the existence of natural a *-homomorphism from $\CovAlg \Ap$ to $\Onmp$:

\state Proposition \label IntroPsi The map $\Phi$ of \lcite{\MapOnTop} vanishes on the redundancy ideal of \/ $\Top \Ap$.
Consequently there exists a *-homomorphism
  $$
  \Psi : \CovAlg \Ap \to \Onmp,
  $$
  such that
  $\Psi\big(\piu(a)\big) = a$, for all $a\in\Ap$, and $\Psi(\Su) = R$.

\Proof
  Let $(a,k)\in \Ap\times\KV$ be a left-redundancy.  Then, taking \lcite{\SimpleRedundancies} into account, for every $b\in
\Ap$, we have that
  $$
  0 = \big(\pio(a)-k\big)\pio(b)\So.
  $$
  Applying $\Phi$ to this leads to
  $$
  0 = \big(a-\Phi(k)\big)bR.
  $$
  In other words, we have that $\big(a-\Phi(k)\big)\Ap R=0$, and hence by \lcite{\TechLemma} we conclude that
  $$
  0 = a-\Phi(k) = \Phi\big(\pio(a)-k\big).
  \eqno{(\dagger)}
  $$
  In the same way we may prove that $(\dagger)$ holds for right-redundancies, hence concluding the proof.  \endProof

Recall from \lcite{\IntroS} that $\S$ is the multiplicative subsemigroup of $\Onm$ generated by all the $s_i$, all the
$t_j$, as well as their adjoints.  In addition to $\S$, we wish to introduce the following subsets of $\S$:

\definition
  \label VariosNotachoes
  \iaitem
  \aitem We shall denote by $\G$ the subset of $\S$ given by
  $$
  \G = \{s_1, \ldots, s_n,\ t_1, \ldots, t_m\}.
  $$
  \aitem We shall denote by $\F$ the subset of $\S$ given by
  $$
  \F = \{s_it_j^*: i\leq n,\ j\leq m\}.
  $$
  \aitem The subsemigroup of $\S$ generated by $\F\cup\F^*$ will be denoted by $\R$.

Observe that, since $\S$ is an inverse semigroup and since the generating set of $\R$ is self-adjoint, one has that $\R$
is itself an inverse semigroup.

\state Proposition \label GenForOnm
  \izitem
  \zitem $\Onmp$ is generated as a C*-algebra by $\F$.
  \zitem $\Ap$ is generated as a C*-algebra by $E(\R)$, the idempotent semi-lattice of\/ $\R$.

\Proof In order to prove (i), let us temporarily denote by $B$ the closed *-subalgebra of $\Onm$ generated by $\F$.  By
\lcite{\sEqualpsq} we have that
  $$
  s_it_j^* = p s_it_j^* p \in \Onmp,
  \subeqmark pstp
  $$
  and hence $B\subseteq\Onmp$.  In order to prove the reverse inclusion it is clearly enough to prove that
  $$
  z:= px_1 \ldots x_rp \in B,
  $$
  whenever $x_k\in\G\cup\G^*$, for every $k\leq r$.
  If $r=0$, that is, if $z=p$, then
  $$
  z = p = \sumn s_is_i^* =
  \sumn s_iqs_i^* =
  \sumn s_it_1^*t_1s_i^* =
  \sumn s_it_1^*(s_it_1^*)^*
  \in B.
  $$

  In case $r>0$, we claim that, unless $z=0$, the $x_k$'s above must:
  \iaitem\aitem start with an element from $\G$,
  \aitem end in an element from $\G^*$, and
  \aitem alternate elements from $\G$ and $\G^*$.

  \medskip In order to prove (a), suppose by contradiction that $x_1\in\G^*$.  Then $x_1^* = px_1^*q$, by \lcite{\sEqualpsq}, so
$px_1 = pqx_1p =0$, and we would have that $z=0$.  A similar reasoning proves (c). As for (b), if two consecutive terms,
say $x_k$ and $x_{k+1}$, both lie in $\G$, then, again by \lcite{\sEqualpsq}, we would have that $x_kx_{k+1} =
px_kqpx_{k+1}q = 0$, and again $z=0$.

This said, we may rewrite $z$ as
  $$
  z=pu_1v_1^*\ldots u_lv_l^*p,
  $$
  where $u_k,v_k\in\G$ (as opposed to $\G\cup\G^*$).  Therefore it suffices to prove that each $u_kv_k^*\in B$.

Since $u_k$ and $v_k$ may be chosen among the $s_i$'s or the $t_j$'s, we are left with the task of proving that
  $$
  s_it_j^*,\ t_js_i^*,\ s_is_j^*,\ t_it_j^*\in B.
  $$
  It is evident that the first two terms above do lie in $B$, while
  $$
  s_is_j^* =
  s_iqs_j^* =
  s_it_1^*t_1s_j^* =
  s_it_1^*(s_jt_1^*)^* \in B.
  $$
  A similar argument proves that $t_it_j^*\in B$.

Focusing now on (ii), let $C$ be the closed *-subalgebra of $\Onm$ generated by $E(\R)$.  Given any element $e \in E(\R)$,
choose $z \in \R$, such that $e = zz^*$.  Then clearly $e \in E(\S)$, and hence $a\in \A$. Since $z\in\R$, equation
\lcite{\pstp} implies that $z = pz$, so $e=pe$, and hence $e\in p\A = \Ap$.  This shows that $C\subseteq\Ap$.

Recall that $\A$ is generated by $E(\S)$, and consequently $\Ap\ (= p\A)$ is generated by $pE(\S)$.  In order to prove
that $\Ap\subseteq C$, it is therefore enough to prove that
  $$
  pe \in C \for e \in E(\S).
  $$
  Write $e=zz^*$, for some $z\in\S$, and further write
  $
  z = x_1 \ldots x_r,
  $
  where $x_k\in\G\cup\G^*$, for every $k\leq r$.  Summarizing, we must prove that
  $$
  f:= p x_1 \ldots x_r x_r^* \ldots x_1^*p \in C,
  $$
  observing that the extra $p$ on the right-hand side above may be added because $E(S)$ is commutative.  Excluding the
trivial case in which $f=0$, we have already seen that (a)--(c) above must hold.  Depending on whether $r$ is even or
odd, we therefore have two alternatives:
  $$
  f = p\ u_1v_1^*\ldots u_{l-1}v_{l-1}^*\ u_lv_l^* \ v_lu_l^*\ v_{l-1}u_{l-1}^*\ldots v_1u_1^*\ p,
  \subeqmark OneAlt
  $$
  or
  $$
  f = p\ u_1v_1^*\ldots u_{l-1}v_{l-1}^*\ u_l\phantom{v_l^* \ v_l}u_l^*\ v_{l-1}u_{l-1}^*\ldots v_1u_1^*\ p,
  \subeqmark OtherAlt
  $$
  where $u_k,v_k\in\G$.
  However \lcite{\OtherAlt} may easily be reduced to \lcite{\OneAlt}, by plugging $v_l=u_l$, so we may assume
\lcite{\OneAlt}.  In order to conclude the proof, it is now enough to show that $u_kv_k^*\in\R$, for all $k$, which we do
by observing that
  $$
  \vcenter{\ialign{
  \hfil $# \ \in \R,$ \pilar{14pt}\cr
  s_it_j^* \cr
  t_js_i^* \cr
  s_is_j^* = s_it_1^*t_1s_j^* \cr
  t_it_j^* = t_is_1^*s_1t_j^* \cr
  }}
  $$
  and the proof is concluded.
  \endProof

\state Proposition \label FirstOccurencerij For every $i\leq n$, and $j\leq m$, let
  $$
  \def\spc{\kern 8pt}
  p_i = s_is^*_i, \spc
  q_j = t_jt^*_j, \spc
  \hp_i = \piu(p_i),\spc
  \hq_j = \piu(q_j), \spc
  \hbox{ and } \spc
  r_{i,j} = \sqrt{nm}\,\hp_i \Su \hq_j.
  $$
  Then
  $
  \Psi(r_{i,j}) = s_it^*_j.
  $
  Consequently $\Psi$ is surjective.

\Proof
  Given $i$ and $j$, we have
  $$
  \Psi(r_{i,j}) =
  \sqrt{nm}\,\psi\big(\piu(p_i)\Su \piu(q_j)\big) =
  \sqrt{nm}\,p_iST^*q_j =
  \su k n \su l m p_is_kt^*_lq_j =
  s_it^*_j.
  $$
  The last sentence in the statement then follows from \lcite{\GenForOnm.i}.  \endProof

\section Mapping $\Onmp$ into $\CovAlg \Ap$

Our main goal is to prove that $\Psi$ is in fact an isomorphism.  In order to accomplish this we will find a representation
of the generators and relations defining $\Onm$ within the algebra of $2\times2$ matrices over $\CovAlg \Ap$, and then we
will employ the universal property of $\Onm$ to construct an inverse for $\Psi$.  Let us begin by proving some useful
algebraic relations.

\state Proposition \label AssortedRels For every $i\leq n$, $j\leq m$, and $f\in\Ap$, one has
  \izitem \def\maior{\pilar{11pt}}
  \zitem $\Caixa {M(q_j)}               = \frac 1m q$,         \zitemmark Mqjfrac      \maior
  \zitem $\caixa {L(p_i)}               = \frac 1n q$,         \zitemmark Kpifrac      \maior
  \zitem $\caixa {\V(q_j)}              = \frac 1m p$,         \zitemmark Vqj          \maior
  \zitem $\caixa {\H(p_i)}              = \frac 1n p$,         \zitemmark Hpi          \maior
  \zitem $\Caixa {p_i\V\big(\H(p_i f)\big)}     = \frac 1n  p_i f$,    \zitemmark piVHpig      \maior
  \zitem $\caixa {q_j\H\big(\V(q_j f)\big)}     = \frac 1m  q_j f$,    \zitemmark qjHVqjg      \maior
  \zitem $\Caixa {\hp_i \Su \Su^*\hp_i}  = \frac 1n \hp_i$,     \zitemmark OneRedun     \maior
  \zitem $\caixa {\hq_j \Su \Su^*\hq_j}  = \frac 1m \hq_j$.     \zitemmark OtherRedun   \maior

\Proof In order to prove \lcite{\Mqjfrac}, we compute:
  $$
  M(q_j) =
  \frac 1m \su km t_k^* q_j t_k =
  \frac 1m \su km t_k^* t_jt_j^* t_k =
  \frac 1m t_j^* t_j =
  \frac 1m q,
  $$
  while \lcite{\Kpifrac} follows similarly.
  As for \lcite{\Vqj}, we have
  $$
  \V(q_j) = \alpha\big(M(q_j)\big) \explain{\Mqjfrac}=
  \frac 1m \alpha(q) =
  \frac 1m \sumn s_iqs^*_i =
  \frac 1m \sumn s_is^*_i =
  \frac 1m p,
  $$
  proving \lcite{\Vqj}, and a similar argument proves \lcite{\Hpi}.  As for \lcite{\piVHpig}, we have
  $$
  \V\(\H(p_i f)\) =
  \alpha M\beta L(p_i f) =
  \alpha L(p_i f).
  $$
  Notice that
  $$
  L(p_i f) =
  \frac 1n \su kn s_k^*p_i f s_k =
  \frac 1n s_i^*s_is_i^*f s_i =
  \frac 1n s_i^*f s_i.
  $$
  So
  $$
  \alpha\big(L(p_i f)\big) =
  \frac 1n \alpha(s_i^*f s_i) =
  \frac 1n \su kn s_ks_i^*f s_is_k^*,
  $$
  and consequently
  $$
  p_i\alpha\big(L(p_i f)\big) =
  \frac 1n \su kn p_is_ks_i^*f s_is_k^* =
  \frac 1n s_is_i^*s_is_i^*f s_is_i^* =
  \frac 1n p_i f,
  $$
  thus proving \lcite{\piVHpig}, while \lcite{\qjHVqjg} follows from a similar argument.

  \def\kp{\uo p}%
  Focusing on \lcite{\OneRedun},
  and letting $\kp_i = \pio(p_i)$,
  we claim that $(\frac 1n p_i,\kp_i \So\So^*\kp_i)$ is a left-redundancy.  To see this, again taking
\lcite{\SimpleRedundancies} into account, pick $f\in\Ap $.  Then
  $$
  \kp_i\So\So^*\kp_i\pio(f)\So =
  \pio(p_i)\So\So^*\pio(p_i f)\So =
  \pio(p_i)\pio\big(\V\(\H(p_i f)\)\big) \So \$=
  \pio\big(p_i\V\(\H(p_i f)\)\big) \So \explain{\piVHpig}=
  \frac 1n \pio(p_i f) \So =
  \pio(\frac 1n p_i)\pio(f)\So,
  $$
  proving the claim, and hence that $\frac 1n \hp_i=\hp_i\Su\Su^*\hp_i$, in $\CovAlg \Ap$.  The last point is proved
similarly.  \endProof

We next present some important algebraic relations among the elements $r_{i,j}$ introduced in
\lcite{\FirstOccurencerij}.

\state Lemma \label algebraicrij For every $i,k\leq n$, and every $j,l\leq m$, one has that
  \Zitem $r_{i,j} r_{k,l}^* = 0$, if $j\neq l$,
  \zitem $r_{i,j}^* r_{k,l} = 0$, if $i\neq k$,
  \zitem $r_{i,j} r_{i,j}^* = \hp_i$,
  \zitem $r_{i,j}^* r_{i,j} = \hq_j$,

\Proof Point (i) follows from the fact that the $\hq_j$ are pairwise orthogonal projections, while (ii) follows from a
similar assertion about the $\hp_i$.  As for (iii), notice that
  $$
  r_{i,j} r_{i,j}^* =
  nm \hp_i \Su \hq_j \Su^*\hp_i =
  nm \hp_i \Su \piu(q_j) \Su^* \hp_i \$=
  nm \hp_i\piu\big(\V(q_j)\big) \Su \Su^*\hp_i \explain{\AssortedRels.\Vqj}=
  n \hp_i\piu(p) \Su \Su^*\hp_i =
  n \hp_i \Su \Su^*\hp_i \explain{\AssortedRels.\OneRedun}=
  \hp_i,
  $$
  proving (iii).
  The proof of (iv) is similar.  \endProof

We will now describe a representation of the generators and relations defining $\Onm$ within the algebra of $2\times2$
matrices over $\CovAlg \Ap$.

\state Proposition For every $i\leq n$, and $j\leq m$, consider the following elements of $M_2\big(\CovAlg \Ap\big)$:
  $$
  \sigma_i=\lmat{r_{i,1}r_{1,1}^*}=\Caixa{r_{i,1}r_{1,1}^*} \otimes e_{2,1}
  $$
  and
  $$
  \tau_j= \lmat{\caixa{\hfil r_{1,j}^*}} = \caixa{\hfil r_{1,j}^*} \otimes e_{2,1},
  $$
  Then, using brackets to denote Boolean value, we have

  \izitem\noindent
  $
  \def\fa{\pilar{12pt},& \hbox{ for all }}
  \def\crr{\hfill\cr}
  \def\itm{\rm(\rzitemplus) & }
  \narrowmatrix{
  \itm \sigma_i^*\sigma_j &=& \bool{i=j}\hp_1 \otimes e_{1,1} \fa i,j\leq n, \crr
  \itm \tau_i^*\tau_j &=& \bool{i=j}\hp_1 \otimes e_{1,1} \fa i,j\leq m, \crr
  \itm \sigma_i\sigma_i^* &=& \hp_i\otimes e_{2,2}             \fa i\leq n,   \crr
  \itm \tau_j\tau_j^* &=& \hq_j \otimes e_{2,2}           \fa j\leq m.   \crr
  }
  $

\Proof \def\a{\noindent(\rzitemplus)} \izitem

\medskip \a
  $$
  \sigma_i^*\sigma_j =
  r_{1,1}r_{i,1}^* r_{j,1}r_{1,1}^* \otimes e_{1,1} \explain {\algebraicrij.ii} =
  \bool{i=j}r_{1,1}r_{i,1}^* r_{i,1}r_{1,1}^* \otimes e_{1,1} \explain {\algebraicrij.iv} =
  \bool{i=j}r_{1,1}\hq_1 r_{1,1}^* \otimes e_{1,1} \$=
  \bool{i=j}r_{1,1}r_{1,1}^* \otimes e_{1,1} \explain {\algebraicrij.iii} =
  \bool{i=j}\hp_1 \otimes e_{1,1}.
  $$
  \a
  $$
  \tau_i^*\tau_j = r_{1,i}r_{1,j}^*\otimes e_{1,1} \explain {\algebraicrij.i} =
  \bool{i=j}r_{1,i}r_{1,i}^*\otimes e_{1,1} \explain {\algebraicrij.iii} =
  \bool{i=j}\hp_1 \otimes e_{1,1}.
  $$
  \a
  $$
  \sigma_i\sigma_i^* =
  r_{i,1}r_{1,1}^* r_{1,1} r_{i,1}^*\otimes e_{2,2} = \explain {\algebraicrij.iv} =
  r_{i,1} \hq_1 r_{i,1}^*\otimes e_{2,2} =
  r_{i,1} r_{i,1}^*\otimes e_{2,2} \explain {\algebraicrij.iii} =
  \hp_i\otimes e_{2,2}.
  $$
  \a
  $$
  \tau_j\tau_j^* =
  r_{1,j}^* r_{1,j} \otimes e_{2,2} \explain {\algebraicrij.iv} =
  \hq_j \otimes e_{2,2}.
  \endProof

As a consequence we see that the $\sigma_i$ and the $\tau_j$ satisfy the relations in \lcite{\DefineOnm}, with the role of $q$
and $p$ being played, respectively, by $\hp_1 \otimes e_{1,1}$, and $\hat p \otimes e_{2,2}$, where
  $$
  \hat p:= \sumn \hp_i = \summ \hq_j = \piu(p).
  $$

We should remark that the validity of the equation ``$p+q=1$'', appearing in \lcite{\DefineOnm}, is guaranteed by the fact
that the $\sigma_i$ and the $\tau_j$ lie in the corner of $M_2\big(\CovAlg \Ap\big)$ determined by the projection
  $$
  \hp_1 \otimes e_{1,1} + \hat p \otimes e_{2,2} = \pmatrix{\hp_1 \otimes e_{1,1} & 0 \cr0 & \hat p \otimes e_{2,2}}.
  $$

  The universal property of $\Onm$ therefore yields:

\state Corollary
  There exists a (not necessarily unital) *-homomorphism
  $$
  \Gamma: \Onm \to M_2\big(\CovAlg \Ap\big)
  $$
  such that
  $$
  \Gamma(s_i) = \sigma_i
  \and
  \Gamma(t_i) = \tau_j,
  \subeqmark ValuesOfGamme
  $$
  for all $i\leq n$, and all $j\leq m$.

Since we are mostly interested in the subalgebra $\Onmp$ of $\Onm$, it is useful to understand the behavior of $\Gamma$ on
this subalgebra.

\state Proposition For all $i\leq n$, and all $j\leq m$, one has that
  $$
  \Gamma(s_it^*_j) = r_{i,j} \otimes e_{2,2}.
  $$
  Consequently the image of $\Onmp$ under $\Gamma$ is contained in the corner of $M_2\big(\CovAlg \Ap)$ determined by $e_{2,2}$.

\Proof
  By \lcite{\ValuesOfGamme}, we have that
  $$
  \Gamma(s_it^*_j) = \sigma_i\tau^*_j =
  r_{i,1}r_{1,1}^* r_{1,j} \otimes e_{2,2}.
  $$
  We must therefore compute
  $$
  r_{i,1}r_{1,1}^* r_{1,j} =
  (nm)^{3/2}\ \hp_i \Su \hq_1 \Su^* \hp_1 \Su \hq_j =
  (nm)^{3/2}\ \piu(p_i) \Su \piu(q_1)\Su^* \piu(p_1)\Su \piu(q_j) \$=
  (nm)^{3/2}\ \piu\big(p_i \V(q_1)\big)\Su \piu\big(\H(p_1)q_j\big) \explain{\AssortedRels.\Vqj\&\Hpi}=
  (nm)^{3/2-1}\piu(p_i)\Su \piu(q_j) =
  \sqrt{nm}\,\hp_i\Su\hq_j = r_{i,j},
  $$
  concluding the calculation of $\Gamma(s_it^*_j)$.  The last assertion in the statement now follows from
\lcite{\GenForOnm.i}.
  \endProof

Observing that the corner of $M_2\big(\CovAlg \Ap)$ determined by $e_{2,2}$ is naturally isomorphic to $\CovAlg \Ap$, we
deduce from the above that:

\state Corollary There exists a *-homomorphism
  $$
  \Lambda: \Onmp \to \CovAlg \Ap,
  $$
  such that $\Lambda(s_it^*_j) = r_{i,j}$.  Moreover,
  $$
  \Gamma(a) = \Lambda(a)\otimes e_{2,2}
  \for a\in\Onmp.
  $$

We are now ready for our main result.

\state Theorem The homomorphism $\Lambda$ of the Corollary above is the inverse of the homomorphisms $\Psi$ of \lcite{\IntroPsi},
and hence $\Onmp$ is *-isomorphic to $\CovAlg \Ap$.

\Proof
  By \lcite{\FirstOccurencerij} we have that $\Psi(r_{i,j}) = s_it^*_j$, and, as seen above, $\Lambda(s_it^*_j) = r_{i,j}$.
Therefore
  $\Psi\circ\Lambda$ acts like the identity on the $s_it^*_j$, and hence
  $$
  \Psi\circ\Lambda = id_{\Onmp},
  \subeqmark PsiLambdaId
  $$
  by \lcite{\GenForOnm.i}.

The proof will then be concluded once we prove that $\Lambda$ is surjective.
  With this goal in mind, we first claim that the $r_{i,j}$ \"{normalize} $\hAp := \piu(\Ap)$, in the sense that
  $$
  r_{i,j}\hAp r_{i,j}^*\subseteq\hAp
  \and
  r_{i,j}^*\hAp r_{i,j}\subseteq\hAp.
  $$
  To see this, let $f\in\Ap$ and observe that
  $$
  r_{i,j}\piu(f)r_{i,j}^* =
  nm\hp_i\Su\hq_j\piu(f)\hq_j\Su^*\hp_i =
  nm\hp_i\piu\big(\V(q_j f q_j)\big)\Su\Su^*\hp_i \$=
  nm\piu\big(\V(q_j f )\big)\hp_i\Su\Su^*\hp_i \explain{\AssortedRels.\OneRedun}=
  m\piu\big(\V(q_j f )\big)\hp_i \in \hAp,
  $$
  and similarly that $r_{i,j}^*\piu( f ) r_{i,j}\in\hAp$.  As a consequence we deduce that any element of $\CovAlg \Ap$,
which is a product of some $r_{i,j}$ and their adjoints, also normalizes $\hAp$.  For any $z$ in the semigroup $\R$
introduced in \lcite{\VariosNotachoes.c}, we have that $\Lambda(z)$ is such a product, so it follows that
  $$
  \Lambda(z)\hAp \Lambda(z)^* \subseteq \hAp,
  $$
  and, in particular, that
  $$
  \Lambda(z)\Lambda(z)^* = \Lambda(zz^*) \in \hAp.
  $$
  Consequently, the idempotent semi-lattice of $\R$ is mapped into $\hAp$ by $\Lambda$.  By \lcite{\GenForOnm.ii}, we then
conclude that
  $$
  \Lambda(\Ap)\subseteq\hAp.
  \subeqmark LambdaApHAp
  $$

By \lcite{\IntroPsi}, we have that $\Psi\big(\piu(a)\big) = a$, for all $a\in\Ap$, and hence $\Psi(\hAp)\subseteq\Ap$.  Suitably restricted, we
may therefore view $\Psi$ and $\Lambda$ as maps
  $$
  \Lambda|_{\Ap}:\Ap \to \hAp
  \and
  \Psi|_{\hAp}:\hAp\to\Ap.
  $$
  By \lcite{\PsiLambdaId} it is clear that
  $$
  \Psi|_{\hAp}\circ\Lambda|_{\Ap} = id_{\Ap}.
  \subeqmark PsiLambdaNoApe
  $$

On the other hand, again by \lcite{\IntroPsi}, we have that
  $$
  \Psi|_{\hAp}\circ\piu = id_{\Ap}.
  \subeqmark Psipiu
  $$
  Viewing $\piu$ as a map
  $$
  \piu: \Ap \to \hAp,
  $$
  notice that \lcite{\Psipiu} implies that $\piu$ is injective, but since it is also clearly surjective, we deduce that
$\piu$ is an isomorphism from $\Ap$ onto $\hAp$.  Once more employing \lcite{\Psipiu}, we conclude that
  $\Psi|_{\hAp}$ is the inverse of $\piu$, and hence it is also an isomorphism.  Using \lcite{\PsiLambdaNoApe} then implies
that $\Lambda|_{\Ap}$ is an isomorphism as well and, in particular, that \lcite{\LambdaApHAp} is in fact an equality of sets.

This said we therefore see that $\hAp$ is contained in the range of $\Lambda$.  Since $\CovAlg \Ap$ is generated by $\hAp$ and
$\Su$, in order to prove our stated goal that $\Lambda$ is surjective, it now suffices to check that $\Su$ lies in the range
of $\Lambda$.  But this follows easily from the fact that $p$ is the unit of $\Ap$ and hence that
  $$
  \Su = \piu(p)\Su\piu(p) =
  \sumn\summ \hp_i \Su \hq_j =
  \frac 1 {\sqrt{nm}}\sumn\summ r_{i,j} =
  \frac 1 {\sqrt{nm}}\sumn\summ \Lambda(s_it^*_j)
  \in \Lambda(\Onmp).
  \endProof

\references

\def \bib #1#2#3#4#5{\bibitem{#1}{#3}{#4}{#5}}

\bib {\AraExel}{AEK}
  {P. Ara, R. Exel and T. Katsura}
  {Dynamical systems of type (m,n) and their C*-algebras}
  {\it Ergodic Theory Dynam. Systems, \rm to appear, arXiv:1109.4093v1}

\bib {\Arz}{A}
  {V.  Arzumanian}
  {Operator algebras associated with polymorphisms}
  {\it Zap. Nauchn. Sem. S.-Peterburg. Otdel. Mat. Inst. Steklov. (POMI) \bf 326 \rm (2005), Teor. Predst. Din. Sist. Komb. i
Algoritm. Metody. 13, 23--27, 279; translation in \it J. Math. Sci.  (N. Y.) \bf 140 \rm (2007), 354--356}

\bib {\ArzRena}{AR}
  {V. Arzumanian and J. Renault}
  {Examples of pseudogroups and their C*-algebras}
  {in {\it Operator algebras and quantum field
   theory (Rome, 1996)}, 93--104, Internat. Press, Cambridge, MA, 1997}

\bib {\RieffelMorita}{BGR}
  {L. G. Brown, P. Green and M. A. Rieffel}
  {Stable isomorphism and strong Morita equivalence of C*-algebras}
  {\it Pacific J. Math., \bf 71 \rm (1977), 349--363}

\bib {\CuVer}{CV}
  {J. Cuntz and A. Vershik}
  {C*-algebras associated with endomorphisms and polymorphisms of compact abelian groups}
  {Preprint, arXiv:1202.5960}

\bib {\tpa}{E1}
  {R. Exel}
  {Twisted partial actions, a classification of regular C*-algebraic bundles}
  {\it Proc. London Math. Soc., \bf 74 \rm (1997), 417--443, arXiv:funct-an/9405001}

\bib {\amena}{E2}
  {R. Exel}
  {Amenability for Fell bundles}
  {\it J. reine angew. Math., \bf 492 \rm (1997), 41--73, arXiv:funct-an/9604009}

\bib {\inverse}{E3}
  {R. Exel}
  {Partial actions of groups and actions of inverse semigroups}
  {\it Proc. Amer. Math. Soc., \bf 126 \rm (1998), 3481--3494, arXiv:funct-an/9511003}

\bib {\ortho}{E4}
  {R. Exel}
  {Partial representations and amenable Fell bundles over free groups}
  {\it Pacific J. Math., \bf 192 \rm (2000), 39--63, arXiv:funct-an/9706001}

\bib {\ExelEndo}{E5}
  {R. Exel}
  {A new look at the crossed-product of a C*-algebra by an endomorphism}
  {\it Ergodic Theory Dynam. Systems, \bf 23 \rm (2003), 1733--1750, arXiv:math.OA/0012084}

\bib {\interactions}{E6}
  {R. Exel}
  {Interactions}
  {\it J. Funct. Analysis, \bf 244 \rm (2007), 26--62, arXiv:math.OA/0409267}

\bib {\newsgrp}{E7}
  {R. Exel}
  {A new look at the crossed-product of a C*-algebra by a semigroup of endomorphisms}
  {\it Ergodic Theory Dynam. Systems, \bf 28 \rm (2008), 749--789, arXiv:math.OA/0511061}

\bib {\oberwolfach}{E8}
  {R. Exel}
  {Endomorphisms, transfer operators, interactions and interaction groups}
  {Mini-Workshop: Endomorphisms, Semigroups and C*-Algebras of Rings, Mathematisches Forschungsinstitut Oberwolfach
Report No. 20/2012, J. Cuntz, W. Szymanski, and J. Zacharias editors, 2012}

\bib {\infinoa}{EL}
  {R. Exel and M. Laca}
  {Cuntz-Krieger algebras for infinite matrices}
  {\it J. reine angew. Math., \bf 512 \rm (1999), 119--172, arXiv:funct-an/9712008}

\bib {\poly}{ER}
  {R. Exel and J. Renault}
  {Semigroups of local homeomorphisms and interaction groups}
  {\it Ergodic Theory Dynam. Systems, \bf 27 \rm (2007), 1737--1771, arXiv:math.OA/0608589}

\bib {\Zettl}{Z}
  {H. Zettl}
  {A characterization of ternary rings of operators}
  {\it Adv. Math., \bf 48 \rm (1983), 117--143}
  \endgroup

\bigskip {\it E-mail address: \tt ruyexel@gmail.com}

\bye

\bye